\documentclass[11pt]{article}

\usepackage[psamsfonts]{amssymb}
\usepackage{amsfonts,euscript}
\usepackage{amsmath,amsthm,graphicx}
\usepackage{color}
\usepackage[colorlinks]{hyperref}

\newtheorem{thm}{\bf Theorem}
\newtheorem{lem}[thm]{\bf Lemma}
\newtheorem{cor}[thm]{\bf Corollary}
\newtheorem{prop}[thm]{\bf Proposition}
\newtheorem{dfn}[thm]{\bf Definition}

\newcommand{\rr}{\mathbb{R}}
\newcommand {\nn} {\mathbb{N}}
\newcommand {\zz} {\mathbb{Z}}
\newcommand {\bbc} {\mathbb{C}}

\newcommand {\al} {\alpha}

\newcommand {\da} {\delta}
\newcommand {\Da} {\Delta}
\newcommand {\ga} {\gamma}

\newcommand{\om}{\omega}

\newcommand {\sa} {\sigma}

\newcommand {\fy} {\varphi}

\newcommand{\ep}{\varepsilon}

\newcommand{\IN}{{\subset}}

\newcommand {\mmm}{{\setminus}}

\newcommand{\8}{{\infty}}
\newcommand{\io}{{I^\infty}}
\newcommand{\ia}{{I^*}}
\newcommand{\0}{{\varnothing}}
\newcommand{\vse}{$\blacksquare$}

\newcommand{\bj}{{\bf {j}}}
\newcommand{\bi}{{\bf {i}}}
\newcommand{\bk}{{\bf {k}}}

\newcommand{\eV}{{\EuScript V}}
\newcommand{\eS}{{\EuScript S}}

\newcommand{\eG}{{\EuScript G}}
\newcommand{\eK}{{\EuScript K}}
\newcommand{\eF}{{\EuScript F}}

\newcommand{\eB}{{\EuScript B}}

\def \sup {\mathop{\rm sup}\nolimits}
\def \fix {\mathop{\rm fix}\nolimits}
\def \Lip {\mathop{\rm Lip}\nolimits}
\def \min {\mathop{\rm min}\nolimits}

\title{Twofold Cantor sets in $\rr$}

\author{ Kirill Kamalutdinov \and 
Andrey Tetenov}

\begin{document}

\maketitle

\bigskip

\begin{abstract} We introduce a class of self-similar sets which we call {\em twofold Cantor sets}  $K_{pq}$ in $\mathbb R$ which  are totally disconnected, do not have weak separation property and at the same time have isomorphic self-similar structures. \end{abstract}

\section{Introduction}
If a self-similar set does not  possess weak separation property it can have unpredicted and surprising properties, especially if it satisfies some additional regularity conditions.
For example, as it was shown in 2006  by one of the authors in  \cite{Tet06, Tetenov}, a self-similar structure $(\gamma,\eS)$  on a Jordan arc  $\gamma$ in $\rr^2$, which does not satisfy WSP, is possible only if $\ga$ is a line segment and two self-similar structures $(\gamma_1,\eS_1)$ and $(\gamma_2,\eS_2)$ on segments $\gamma_1$ and  $\gamma_2$, which do not satisfy WSP, are isomorphic if and only if  the homeomorphism $\varphi:\ \gamma_1\to\gamma_2$, which induces the the isomorphism of these structures, is a linear map. So the question arises, does such rigidity phenomenon occur for self-similar sets whose dimension is smaller than 1 and which are therefore totally disconnected?\\

 We introduce   new objects which we call {\em twofold Cantor sets}  $K_{pq}$ in $\mathbb R$
 %(or $\mathbb C$)
, which are generated by the systems $\eS_{pq}=\{S_1,S_2,S_3,S_4\}$ consisting of two symmetric pairs of similarities  $S_1(x)=px$, $S_2(x)=qx$ and $S_3(x)=px-p+1$, $S_4(x)=qx+q-1$, with fixed points at 0 and 1 and real (or complex) coefficients p and q such that $\dim_H(K_{pq})<1/2$ and the sets $S_1^{m_1}S_2^{n_1}(A)$ and $S_1^{m_2}S_2^{n_2}(A)$, where $A=S_3(K_{pq})\cup S_4(K_{pq})$,
 are disjoint if $(m_1,n_1)\neq (m_2,n_2)$.\\

Currently, we confine ourselves to the case when $p,q$ are real  and lie in the set $\eV=(0,1/16)^2$.

We show that all such systems $\eS_{pq}$ do not possess weak separation property (Proposition \ref{nowsp}) and  any two of them define isomorphic self-similar structures $(K_{pq},\eS_{pq})$, $(K_{p'q'},\eS_{p'q'})$ (Theorem \ref{rigid}).   At the same time, each homeomorphism $f: K_{pq}\to K_{p'q'}$, which induces the isomorphism of respective self-similar structures, cannot be extended continuously to a homeomorphism   of $\rr$ to itself.  Nevertheless, each $f$ extends to a homeo\-mor\-phism $\tilde f: \bbc\to\bbc$, but such homeomorphism $\tilde f$ never agrees with the systems $\eS_{pq},\eS_{p'q'}$ if $(p,q)\neq(p',q')$ (Theorem \ref{homC}).\\

Though twofold Cantor sets $K_{pq}$ are homeomorphic to the Cantor set, they densify near $0$ and in the case of positive  $p,q$ there is a topological limit $\lim\limits_{t\to+\8}tK_{pq}$, which is  $[0,+\8)$ (Theorem \ref{limit}).\\

All the mentioned properties of twofold Cantor sets follow naturally from their definition. There is the most difficult question   left: do such sets exist? \\

The answer to this question requires different techniques.
 Our problem is how to find those $p,q$, for which each of the intersections $S_1^m(A)\cap S_2^n(A)$ is empty, so we analyse how large is  the set of those pairs $(p,q)$ which do not possess such property. First we show that for any pair of non-negative $m,n\in\zz$ and for any $p\in (0,1/16)$ the set $\Delta_{mn}(p)$ of those $q\in (0,1/16)$, for which   $S_1^m(A)\cap S_2^n(A)\neq\0$, has dimension less than 1 (Theorem \ref{tech3}). Therefore for any  $p$ the set  $\bigcup\limits_{m,n=0}^\8 \Da_{mn}(p)$ has 1-dimensional Lebesgue measure zero in $\{p\}\times(0,1/16)$.\\

To prove that $\dim_H \Delta_{mn}(p)<1$, we use a bunch of two statements, General Position Theorem \ref{genpos} and Displacement Theorem \ref{collage} initially used by the authors in \cite{TKV}. \\

First theorem considers the set of pairs $(\fy_1(\xi,x),\fy_2(\xi,x))$ of $\al$-H\"older maps   $\fy_i:D \times K \to \rr^n$ of a metric compact $K$, depending on a parameter $\xi\in D\IN\rr^n$, and find conditions  under which the set $\{\xi\in D:\fy_1(\xi,K)\cap\fy_2(\xi,K)\neq\0\}$  has  dimension smaller than $n$. Such conditions are:
\\(1) Both $\fy_i$ are  H\"older continuous    with respect to $x\in K$ and \\(2)  The difference $\fy_1(\xi,x_1)-\fy_2(\xi,x_2)$ is inverse Lipschitz with respect to $\xi\in D$.

Lemma \ref{tech2} allows to check the conditions of the Theorem \ref{genpos} for the particular case, when we try to find those values of parameter $q$, for which $S_1^m(A)\cap S_2^n(A)\neq\0$.

This allows us to show that the set $\eK$ of those $(p,q)\in \eV$ which correspond to twofold Cantor sets, has a full measure in $\eV$ (Theorem \ref{fullmes}), while the complement of $\eK$ in $\eV$ is noncountable and dense in $\eV$. \\

Imposing additional non-overlap condition, we can obtain twofold Cantor sets  from the systems  $$\eS=\eF \cup \eG, \ \ \eF=\{F_1,\dots, F_m\},\ \  \eG=\{G_j:\ j\in J\},\ \  J\subseteq \{1,\dots,m\},\fix F_i=\fix G_i\mbox{ for any }i\in J$$  with non-distinct fixed points  studied by Barany \cite{Bar2}. Using   
transversality condition introduced by Simon, Solomyak and Urbanski \cite{Sim1, Sim2} he showed that if $\eF$ satisfies OSC, then for Lebesgue-almost all positive vectors of parameters $(\Lip G_j)_{j\in J}$ in some neighbourhood of zero the Hausdorff dimension $d$ of the attractor $K$  of the system $\eS$ is a solution of  equation $\sum\limits_{i=1}^m (\Lip F_i)^d  +\sum\limits_{j\in J} (\Lip G_j)^d - \sum\limits_{j\in J} (\Lip F_j \Lip G_j)^d  = 1$, and that in such case $H^d(K)=0$.

We prove similar result (Theorem \ref{hdim}) for twofold Cantor sets  using much more simple approach.

%\blue{There are many separation conditions widely used in fractal geometry. }%для того, чтобы выделить семейства IFS, размерность которых поддается вычислению тем или иным образом. 
%\blue{Open set condition (OSC) was first introduced by Moran \cite{Mor}, and then studied by Hutchinson \cite{Hut}, that is for the system of contraction similarities $\mathcal S=\{S_1,\dots,S_m\}$ in $\mathbb{R}^n$ exists an open set $O$, such that $S_i(O)\subset O$ and $S_i(O)\cap S_j(O)=\varnothing$ for all $i,j\in\{1,\dots,m\}$. Moran showed, that if OSC holds, then the Hausdorff dimension of attractor of the system $\mathcal S$ coincides with the similarity dimension of $\mathcal S$, that is the unique positive solution $d$ of Moran equation $\sum\limits_{i=1}^m (\Lip S_i )^d=1$. The relations between OSC, positive Hausdorff measure and other properties were studied by Schief \cite{Schief}, Bandt and Graf \cite{SSS7}}
%
%
%\blue{ }
%
%
%
%
%
%

\section{Twofold Cantor sets}

\subsection{Definition and basic properties of twofold Cantor sets}
\subsubsection{\bf Self-similar sets.}

Let $(X, d)$ be a complete metric space. 
A mapping $S: X \to X$ is a contraction if $\Lip S < 1$ and it  is called a similarity if $d(S(x), S(y)) = r d(x, y)$ for all $x, y\in X$ and some fixed r. 

\begin{dfn} 
Let $\eS=\{S_1, S_2, \ldots, S_m\}$ be a system of  contractions of a complete metric space $(X, d)$.
 A nonempty compact set $K\IN X$ is called the attractor of the system $\eS$, if $K = \bigcup \limits_{i = 1}^m S_i (K)$. \end{dfn}
By Hutchinson's Theorem \cite{Hut}, the attractor $K$ is uniquely defined by the system $\eS$. 
 
 We also call the   set $K$ {\em self-similar} with respect to $\eS$ and the pair $(K,\eS)$ is called a self-similar structure.
Throughout  this paper, the maps $S_i\in \eS$ will be  similarities and the space $X$ will be $\rr$.\smallskip

Given a system $\eS=\{S_1,...,S_m\}$, $I=\{1,2,...,m\}$ is the set of indices, $\ia=\bigcup\limits_{n=1}^\8 I^n$ 
is the set of all finite $I$-tuples, or multiindices $\bj=j_1j_2...j_n$. By $\bi\bj$ we denote the concatenation of the corresponding multiindices; we write $\bi\sqsubset\bj$, if  $\bj=\bi\bk$ for some $\bk\in\ia$.

We write
$S_\bj=S_{j_1j_2...j_n}=S_{j_1}S_{j_2}...S_{j_n}$ 
and for the set $A
\subset X$ we denote $S_\bj(A)$ by $A_\bj$; 
$I^{\8}=\{{\bf \al}=\al_1\al_2\ldots:\ \ \al_i\in I\}$ is the
index space;  $\pi:I^{\8}\rightarrow K$ is the {\em index map}, which sends $\bf\al\in I^{\8}$ to  the point $\bigcap\limits_{n=1}^\8 K_{\al_1\ldots\al_n}$.

We denote by $F_\eS=\{S_\bj: \bj\in\ia\}$ (or simply by $F$)  the semigroup, generated by $\eS$; then $\eF=F^{-1}\circ F$, or a set of all compositions $S_\bj^{-1}S_\bi$ is the associated family of similarities, first considered  by Bandt and Graf \cite{SSS7}. 
Weak separation property (WSP) was introduced by Lau and Ngai \cite{Lau}. It was shown by Zerner \cite{Zer}, that it is equivalent to the topological condition: system $\mathcal S=\{S_1,\dots, S_m\}$ of contraction similarities satisfy the WSP iff $\rm Id \notin \overline{\eF \setminus \rm Id}$.\\

 Notice that throughout this paper $m=4$, therefore $I=\{1,2,3,4\}$.

\medskip

\subsubsection{\bf Systems $\eS_{pq}$ and their attractors.}

 We define a system $\eS_{pq}=\{S_1,S_2,S_3,S_4\}$ of contraction similarities on $[0,1]$ by the equations 
$S_1 (x)=px$, $S_2 (x)=qx$, $S_3 (x)=px+1-p$, $S_4 (x)=qx+1-q$, where the contraction ratios $p,q\in (0,1/2)$.\\
Let  $K_{pq}$ be the attractor of the system $\eS_{pq}$. We write $K$ instead of $K_{pq}$ if it does not cause ambiguity. We will also denote $A=S_3(K) \cup S_4(K)$ and $B=S_1(K) \cup S_2(K)$.\\
Notice that  if $q=p^n,n\in \nn$, the attractor $K_{pq}$ is the Cantor  set $K_p$ generated by $S_1$ and $S_3$.

\medskip

There are several obvious properties of systems $\eS_{pq}$ and their attractors $K_{pq}$ :

\begin{prop}\label{props}\ \\
(i) $S_1 S_2 = S_2 S_1$, $S_3 S_4 = S_4 S_3$;\\
(ii) $K_{pq}=\om(K_{pq})$, where $\om(x)=1-x$;\\
(iii) For  $i=1,2$ and for any integers $m\neq n$,  $S_i ^m (A)\cap S_i ^n (A)=\0$;\\ 
(iv) For any $m,n\in\mathbb N\cup\{0\}$,
$S_1 ^m S_2 ^n (K_{pq})  \subseteq S_1 ^m (K_{pq})\cap S_2 ^n (K_{pq})$; \\
(v) $K_{pq}\setminus\{0\}=\bigcup\limits_{m,n=0}^\infty S_1 ^m S_2 ^n(A)$.\\  
\end{prop}
{\bf Proof:  }\ 
Commutativity (i) and symmetry with respect to $1/2$ (ii) are obvious;   notice, that for any $i\in I$, $\om\circ S_i \circ \om =S_{i'}$ where $i'\equiv i+2(\mod 4)$.\\ 
(iii)   Suppose $m<n$.  $S_i^m(A)\IN (S_i^m(1/2),S_i^m(1)]$, and $S_i^n(A)\IN (S_i^n(1/2),S_i^n(1)]$. Since $S_i^n(1)\le S_i^m(1/2)$, $S_i^m(A)\cap S_i^n(A)=\0$.\\
(iv) is trivial. \\
(v) Let $\sa=i_1i_2...\in \io$ and $x=\pi(\sa)$. Notice that if $\min\{l:i_l\in \{3,4\}\}= k$, then $x\in S_1 ^m S_2 ^n(A)$ for some $m,n\ge 0$ such that  $m+n=k-1$. If $x\notin \bigcup\limits_{m,n=0}^\infty S_1 ^m S_2 ^n(A)$, then  $\sa\in \{1,2\}^\8$, therefore $x=0$.
%%%There was other proof of (v), but not clear:
%The set $\bigcup\limits_{m,n=0}^\infty S_1^mS_2^n(A)$ is dense in $K$. 
%The family of compact sets $\{S_ 1^m S_2^n(A)\}$ is locally finite in $(0,1]$, 
%therefore its complement in $K$ is $\{0\}$.  
\vse\\

{\bf Notation:} We use $\bigsqcup$ to denote disjoint union.

% Twofold Cantor sets
\begin{prop}\label{disjoint}
Let the system $\eS_{pq}$ satisfy the condition

{\bf (TF):} \quad For any $m,n\in\mathbb N$, $S_1 ^m(A)\cap S_2 ^n(A)=\0$. \quad
Then

(i) $K=\{0\}\cup\bigsqcup\limits_{m,n=0}^{\infty}S_1^m S_2^n (A)$;

(ii) for any $m,n\in\mathbb N\cup\{0\}$, $S_1^m (K)\cap S_2^n (K)=S_1^m S_2^n (K)$ .
\end{prop}
{\bf Proof:  }
Notice that for any integers $k,l,m,n$, $$S_1^{m+k}S_2^n(A)\cap S_1^{m}S_2^{n+l}(A)\neq\0\mbox{\quad iff \quad }S_1^k(A)\cap S_2^l(A)\neq\0.$$ 

Therefore {\bf (TF)}  is equivalent to (i). To prove (ii), notice  that if all the sets $S_1^m S_2^n(A)$ are disjoint, then the set $(S_1 ^m (K)\cap S_2 ^n (K))\setminus\{0\}$ is equal to \begin{equation*}\bigcup\limits_{k,l=0}^\infty S_1 ^{m+k} S_2 ^l(A) \cap \bigcup\limits_{k,l=0}^\infty S_1 ^{k} S_2 ^{n+l}(A)=\bigcup\limits_{k,l=0}^\infty S_1 ^{m+k} S_2 ^{n+l}(A)=S_1^m S_2^n (K)\setminus\{0\}.\qquad \blacksquare\end{equation*}

\begin{dfn}\label{deftf}
If the system $\eS_{pq}$  satisfies the condition {\bf TF}, we call $K_{pq}$ a {\em twofold Cantor set}.
\end{dfn}  
\medskip

\subsubsection{\bf Dimension theorem for twofold Cantor sets.}

\begin{thm}\label{hdim}  Hausdorff dimension $\dim_H K_{pq}$ of a twofold Cantor set $K_{pq}$ satisfies the equation \qquad $p^d+ q^d- (pq)^d=1/2$.
\end{thm}
{\bf Proof:  }\
Consider the set $B=S_1(K)\cup S_2(K)$.
By Proposition \ref{props} (ii), $A=\om(B)$, so
from Proposition \ref{disjoint} we get \\
\begin{equation}\label{bunion }B=\{0\}\cup\bigcup\limits_{m+n>0}S_1^m S_2^n A = S_1 (B) \cup S_1 \om(B) \cup \bigcup\limits_{n=1}^{\infty} S_2^n \om(B)\end{equation}

Thus  $B$ is the attractor of infinite system  $\eB=\{S_1, S_1 \om, S_2 \om, S_2^2 \om, \dots, S_2^n \om, \dots \}$ of similarities whose contraction ratios are $\{p,p,q,q^2,\dots,q^n,\dots\}$ respectively. 

Then standard argument \cite[Theorem 9.3]{Fal} shows that  if $d$ is a solution of the equation
 \begin{equation}\label{sereqn}2p^d+\sum\limits_{k=1}^\8q^{nd}=1,\end{equation} then $\dim_HB\le d$. 
 
 In our case
the equation (\ref{sereqn}) is equivalent to the equation $p^d+q^d-(pq)^d=1/2$. If $p,q\in(0,1)$  it has unique positive solution which we denote by $d_{pq}$.

Let $B_n$ be the attractor of a  subsystem $\eB_n =\{S_1, S_1 \om, S_2 \om, \dots, S_2^n \om\}$ of the system $\eB$.

For any $n\in \nn$, $B_n \subset B_{n+1}$ and $\bigcup\limits_{n=1}^{\infty} B_n\subseteq B$, therefore $\dim_H B_n < \dim_H B\le d_{pq}$.
 
Since $K_{pq}$ is a twofold Cantor set, the compact sets $S_1 (B_n)$, $S_1 \om(B_n)$ и $S_2^k \om(B_n)$, $k=1,\dots,n$ are disjoint. Therefore the system $\eB_n$ satisfies open set condition with open set $O=V_\ep(B_n)$, where $ \ep$ is less than half minimal distance between the copies of $B_n$.

 Then Hausdorff dimension of the set $B_n$ is the unique positive solution $d_n$ of the equation $2p^{x}+ \sum\limits_{k=1}^{n} q^{kx}=1$ and the set $B_n$ has positive finite measure in the dimension $d_n$.\\

The sequence $d_n$ increases and $d_n<d_{pq}$, so it has a   limit which satisfies the equation
$p^d+ q^d- (pq)^d=1/2$. Therefore $d_{pq}=\dim_H B=\dim_H K_{pq}$. \vse

\subsection{Density properties and violation of WSP} 
\subsubsection{\bf Logarithmic incommensurability of p and q}
For twofold Cantor sets $K_{pq}$ the logarithms of $p$ and $q$
are incommensurable, which causes most of their unusual properties.\\

\begin{prop}\label{log} 
If $K_{pq}$ is a twofold Cantor set, then $\dfrac{\log p}{\log q}\notin\mathbb Q$.
\end{prop}
{\bf Proof:  }
Assume the contrary. Then for some $m,n$, $p^m=q^n$ implies $S_1^m ({1})= S_2^n (1)$, which contradicts the condition {\bf(TF)}.
\vse

All the statements of this subsection require  the condition $\dfrac{\log p}{\log q}\notin\mathbb Q$ only. Therefore they are valid for twofold Cantor sets.

\medskip

\begin{lem}[see \cite{TKV0}]\label{dense}
If $\dfrac{\log p}{\log q}\notin\mathbb Q$, then the multiplicative  group $G$, generated by $p$ and $q$ is a  dense subgroup in $\langle \mathbb R_{+}, \cdot \rangle$. \vse
\end{lem}

If we expand a twofold set $K_{pq}$, multiplying it by $t$, tending to $+\8$, then a topological limit \cite{Kur} of the sets $tK_{pq}=\{tx, x\in K_{pq}\}$ will be the set $[0,+\8)$:\\   
{\bf Notation.} Let $d_H(X_1,X_2)$ be  Hausdorff distance between non-empty compacts $X_1,X_2\subset \mathbb R$. \quad Denote $\Da_{[a,b]} (X)=d_H(X\cap[a,b],[a,b])$. 

\medskip

\begin{thm}\label{limit}
If $\dfrac{\log p}{\log q}\notin\mathbb Q$, then there is a topological limit $\lim\limits_{t\to +\infty}t K_{pq}=[0,+\8)$.
\end{thm}

{\bf Proof:  }
First we prove that $\lim\limits_{t\to +\infty} \Da_{[0,r]}(t K_{pq})=0$ for any $r>0$.

Put $G=\{p^m q^n:\ m,n\in \mathbb Z\}$ and $G_{+} = \{p^m q^n:\ m,n\in \mathbb N\}$. Take some $r>0$. \\
Consider the sets $G_k = (pq)^{-k} G_{+}$. For any $k\in\nn$, $G_k\subset G_{k+1}$ and $G=\bigcup \limits_{k=0}^{+\infty} G_k$. By Lemma \ref{dense} the set $G$ is dense in $[0,+\infty)$, therefore $\Da_{[0,r]}(G_k) \to 0$ as $k\to\infty$. \\
For any $t>1$ there is an integer $k\ge 0$, such that $(pq)^{-k}\le t < (pq)^{-k-1}$. Notice that $\Da_{[0,r]}(t G_{+})\le \Da_{[0,tr]}(t G_{+})=t \Da_{[0,r]}(G_{+})$, therefore $\Da_{[0,r]}(t G_{+})\leq\dfrac{1}{pq} \Da_{[0,r]}(G_k)$.
Then $\Da_{[0,r]}(t K_{pq})\to 0$ while $t\to +\infty$.

This means that for any $x\in [0,+\8)$ and any $\ep>0$ there is such $M>0$, that if $t>M$, then $d(x,t K_{pq})<\ep$. Therefore $x\in\lim\limits_{t\to +\infty}t K_{pq}$.
\vse\\

The same property holds for all points $S_\bj(0)$ and $S_\bj(1)$:\\

\begin{cor}
If $\dfrac{\log p}{\log q}\notin\mathbb Q$, then for any $r>0$ and any ${\bf i} \in I^{*}$  $$\lim\limits_{t\to +\infty} \Da_{[0,r]}(t (K_{pq}-S_{\bf i}(0)))=\lim\limits_{t\to +\infty}  \Da_{[-r,0]}(t (K_{pq}-S_{\bf i}(1)))=0. \qquad\blacksquare$$
\end{cor}\ \\

The logarithmic incommensurability of $p$ and $q$ causes the violation of WSP for the system $\eS_{pq}$:\\

\begin{prop}\label{nowsp} 
If $\dfrac{\log p}{\log q}\notin\mathbb Q$, then the system $\eS_{pq}$ does not have weak separation property (WSP). 
\end{prop}
{\bf Proof:  }
By Lemma \ref{dense} there is a sequence of $(m_k,n_k)\in \mathbb N^2$, such that $\dfrac{p^{m_k}}{q^{n_k}} \rightarrow 1$ as $k\to\infty$, therefore $S_1 ^{m_k} (S_2 ^{n_k})^{-1} \rightarrow \rm Id$. Since $\dfrac{p^{m_k}}{q^{n_k}}\neq 1$, the point $\rm Id$ is a limit point of $F^{-1} F$. 
\vse \ \\

\bigskip

\subsection{Isomorphisms of twofold Cantor sets.}

Let  $(K,\{S_1,\dots,S_m\})$, $(K',\{S'_1,\dots,S'_m\})$ be self-similar structures. We say that 
a homeomorphism $f:\ K_{pq}\to K_{p'q'}$,  realises the isomorphism of self-similar structures $(K_{pq},\eS_{pq})$ and $(K_{p'q'},\eS_{p'q'})$, if $f(S_i(x))=S_i'(f(x))$ for any $x\in K$ and any $i\in \{1,...,m\}$. 

% rigidity for Jordan arcs
It was proved \cite{Tetenov} that if the systems $\eS$, $\eS'$ of contraction similarities in $\mathbb C$ with Jordan attractors $K$ and $K'$ do not have WSP, and $f:\ K\to K'$  realises the isomorphism of self-similar structures $(K, \eS)$ and $(K', \eS')$, then $f$ is a linear map, and $K$, $K'$ are straight line segments.\\

The question, to what extent such rigidity phenomenon is valid for self-similar structures whose Hausdorff dimension is smaller than 1, remain open still.
Further we'll try to establish isomorphism of self-similar structures on twofold Cantor sets. First we construct a representation 
of points $x\in K_{pq}$, convenient for establishing  such isomorphisms.
\medskip

{\bf Notation.}  Let $(m,n), (m', n')\in\zz\times\zz$. We will write $(m,n)<(m', n')$ if $m\le m'$, $n\le n'$ and $m+n<m'+n'$.\\

Let $G_+=\{S_1^mS_2^n:(m,n)>(0,0)\}$ and
$H_+=\{S_3^mS_4^n:(m,n)>(0,0)\}$ be  semigroups generated by $S_1,S_2$  and $S_3,S_4$ respectively. Then it follows from  Propositions \ref{props} and \ref{disjoint}  that
\begin{equation} B=\bigsqcup\limits_{g\in G_+}g(A)\qquad A=\bigsqcup\limits_{h\in H_+}h(B) \label{djun1}\end{equation}

The following Lemma represents each element of $K_{pq}$ as a unique alternating sum.\\ 

\begin{lem}\label{series}
Let $K_{pq}$ be a twofold Cantor set. Then each $x\in K_{pq}\mmm\{0\}$ has unique alternating sum representation
\begin{equation}\label{finfsum}x=\sum\limits_{k=0}^{N} (-1)^k p^{m_k} q^{n_k}
\mbox{ \qquad  or \qquad  }x=\sum\limits_{k=0}^{\8} (-1)^k p^{m_k} q^{n_k} \end{equation}
where $(m_k,n_k)<(m_{k+1},n_{k+1})$ for any $k$ and $(m_0,n_0)\ge(0,0)$.\\ 
This sum is finite if $x\in \bigcup\limits_{{\bf i}\in I^{*}}S_{\bf i}(\{0,1\})$ and infinite otherwise.\\
\end{lem}{\bf \noindent Proof:  } Suppose $x=x_0\in B\mmm\{0\}$, then there is unique  $g_0\in G^+$ and  $y_0\in A$, that $x=g_0(y)$.\\
If $y=1$, then $x=g_0(1)$, otherwise there is unique  $h_0\in H^+$ and   $x_1\in B$, such that
$y=h_0(x_1)$.\\ If $x_1=0$, then $x=g_0h_0(0)$, otherwise $x=g_0h_0(x_1)$. Proceeding by induction, we see that for any $x\in B\mmm\{0\}$ there is unique representation having  one of the forms

\begin{equation}\label{3forms}x=g_0h_0...g_Nh_N(0) \mbox{\quad or \quad } x=g_0h_0...g_N(1)
\mbox{ \quad  or \quad  }\{x\}=\bigcap\limits_{N=0}^{\8} g_0h_0...g_Nh_N(K) \end{equation}

 We obtain respective representations for $x\in A$, if we take $g_0={\rm Id}$.

Notice that $gh(x)=S_1^{i} S_2^{j} S_3^{i'} S_4^{j'} (x)=p^i q^j -p^{i+i'} q^{j+j'}+p^{i+i'} q^{j+j'} x$, therefore 
\begin{equation}\label{frep}g_0h_0...g_N(x)=p^{m_{2N}} q^{n_{2N}}x+\sum\limits_{k=0}^{N-1} (p^{m_{2k}} q^{n_{2k}}-p^{m_{2k+1}} q^{n_{2k+1}}),\end{equation}
\begin{equation}\label{infrep}g_0h_0...g_Nh_N(x)=p^{m_{2N+1}} q^{n_{2N+1}}x+\sum\limits_{k=0}^{N} (p^{m_{2k}} q^{n_{2k}}-p^{m_{2k+1}} q^{n_{2k+1}}),\end{equation} where  $m_{2k+1}=\sum\limits_{l=1}^k (i_l+i'_l)$, $n_{2k+1}=\sum\limits_{l=1}^k (j_l+j'_l)$, $m_{2k}=m_{2k-1}+i_{k}$ and $n_{2k}=n_{2k-1}+j_{k}$. \\

Applying formulas (\ref{frep}),(\ref{infrep}) to each of the relations in (\ref{3forms}), we get the desired result (\ref{finfsum}).
\vse

\begin{thm}\label{rigid}
Let $K_{pq}$, $K_{p'q'}$ be twofold Cantor sets. Then:\\
(i) There is a homeomorphism $f:\ K_{pq}\to K_{p'q'}$, which realises the isomorphism of self-similar structures $(K_{pq},\eS_{pq})$ and $(K_{p'q'},\eS_{p'q'})$.\\
(ii)  If $(p,q)\neq(p',q')$, then  $f$ cannot be extended to a homeomorphism of $[0,1]$ to  itself.

\end{thm}
{\bf Proof:  }\ 
(i) It follows from Lemma \ref{series} that each element $x\in K_{pq}$ has unique representation either as a finite sum 
$x=\sum\limits_{k=0}^{N} (-1)^k p^{m_k} q^{n_k}$ or as sum of a series  $\sum\limits_{k=0}^{+\infty} (-1)^k p^{m_k} q^{n_k}$.
The same is true for the set $K_{p'q'}$. Therefore,  one can define a bijection $f:\ K_{pq}\to K_{p'q'}$ by

\[f(x) =
\begin{cases}
\sum\limits_{k=0}^{N} (-1)^k {p'}^{m_k} {q'}^{n_k}, & \text{if  \quad}  x=\sum\limits_{k=0}^{N} (-1)^k p^{m_k} q^{n_k} \\
\sum\limits_{k=0}^{+\infty} (-1)^k {p'}^{m_k} {q'}^{n_k}, & \text{if  \quad}   x=\sum\limits_{k=0}^{+\infty} (-1)^k p^{m_k} q^{n_k}.
\end{cases}
\]

The sets $A$ and $B$ are open-closed in $K_{pq}$. Applying formulas (\ref{djun1}) stepwise we conclude that each of the sets $g_0h_0...g_Nh_N(B)$, $g_0h_0...g_N(A)$ is open-closed in 
$K_{pq}$. 

The sets $U_N=\bigcup\limits_{k=0}^N S_1^kS_2^{N-k}(B)$ are also open-closed, because their complement is a finite union 
$\bigsqcup\limits_{k+l\le N}S_1^kS_2^l(A)$. These sets form a fundamental system of open-closed neighborhoods of 0 in $K_{pq}$. 

Therefore, an open-closed neighborhood base of  every point\\
$x=g_0h_0...g_Nh_N(0)$\quad or \quad $x=g_0h_0...g_N(1)$
is $\{g_0h_0...g_Nh_N(U_k),k\in \nn\}$\quad or \quad$ \{g_0h_0...g_N\om(U_k),k\in \nn\}$ respectively.

For the points having representation
 $x=\bigcap\limits_{N=0}^{\8} g_0h_0...g_Nh_N(K)$, such open-closed neigh\-bour\-hood base is $\left\{ g_0h_0...g_kh_k(B),k\in\nn\right\}$.
 
 The map $f$ sends all these base sets to respective base sets in $K_{p'q'}$. Therefore $f: K_{pq}\to K_{p'q'}$ is a homeomorphism. 
 
  \medskip
(ii) Suppose that the map $f:K_{pq}\to K_{p'q'}$ which induces the isomorphism of self-similar structures $(K_{p q},S_{p q})$ and $(K_{p' q'},S_{p' q'})$ can be extended to a homeomorphism $\tilde f$: $[0,1]\to [0,1]$. Then
 the function $\tilde f$ is increasing. \\
Put $\alpha=\dfrac{\log p'}{\log p}$, $\beta =\dfrac{\log q'}{\log q}$. There are 2 cases: $\alpha \neq \beta$ and $\alpha=\beta$. \\
In the first case suppose $\alpha<\beta$. There is $\lambda\in\mathbb Q$, such that $\alpha<\lambda<\beta$. Take  $k,l,m,n\in \mathbb N$ such that $\dfrac{k-m}{n-l}=\lambda$. Simple computation shows that $p^k q^l < p^m q^n$ and $p'^k q'^l > p'^m q'^n$, which violates monotonicity of $\tilde f$.\\

Let now $\alpha=\beta$. Then for any $x\in G_+$, $ f(x)=x^\al$.\\ 
For any $m,n\in\mathbb N$ the inequalities  $p^m < q^n (1-p)$ and $p'^m < q'^n (1-p')$ are equivalent because $f(p^m)=p'^m$, $q'^n (1-p')=f(q^n (1-p))=S_2 ^n S_3 f(0)$, and $f$ increases. \\

Consider the set $W=\left\{\dfrac{p^m}{q^n}:\ p^m < (1-p)q^n;\ m,n\in\mathbb N\right\}$. The set $W$ is dense in $(0,1-p)$,  so $\sup W= 1-p$. At the same time,\\   $$W'=\left\{\dfrac{p'^m}{q'^n}:\ p'^m < (1-p') q'^n;\ m,n\in\mathbb N\right\}=\{x^\al, x\in W\}$$ Therefore $\sup A'=1-p'=(\sup A)^\al$
or $1-p^\alpha =(1-p)^\alpha$, which implies $\alpha=1$ and $(p, q)=(p', q')$.
\vse

\medskip

\begin{thm}\label{homC}
Let $f:K_{pq}\to K_{p'q'}$ be a homeomorphism of twofold Cantor sets which induces the isomorphism of self-similar structures $(K_{pq},\eS_{pq})$ and $(K_{p'q'},\eS_{p'q'})$. Then $f$ has an extension to a homeomorphism $\tilde f:\bbc\to\bbc$. If $\tilde f\circ S_i(z) = S'_i \circ \tilde f(z)$ for any $z\in\mathbb C$ and for any $i\in I$ then $(p,q)=(p',q')$.
\end{thm}
{\bf Proof:  }\ 
A homeomorphism of  totally disconnected compact sets in $\mathbb C$ has an extension to a self-homeomorphism of $\mathbb C$   \cite[Ch.13, Theorem 7]{Moise}. Let $\tilde f:\ \mathbb C \to \mathbb C$ be the extension of the map $f$ which induces the isomorphism of $(K_{pq},\eS_{pq})$ and $(K_{p'q'},\eS_{p'q'})$. \\Suppose $\tilde f\circ S_i(z) = S'_i \circ \tilde f(z)$ for any $z$ and $i$.  \\
Then for any $m,n\in\zz$, $f(p^m q^n)=p'^m q'^n$, so $f$ sends $G=\{p^m q^n:\ m,n\in\zz \}$ to $G'=\{p'^m q'^n:\ m,n\in\zz \}$. Since $K_{pq}$ and $K_{p'q'}$ are twofold Cantor sets, the groups $G $ and $G'$  are dense in $[0,1]$,  so $f([0,1])=[0,1]$ which contradicts Theorem \ref{rigid}(ii).
\vse

%%%%%%%%%%%%%%%%%%%%%%
\subsection{Existence of twofold Cantor sets}

\subsubsection{\bf General position Theorem.}

Let $K$ be the attractor of a system $\eS=\{S_1,...S_m\}$ of contraction similarities of $\mathbb R^n$, and let $\dim_H K<n/2$. Suppose for some $k,l$,  $S_k(K)\cap S_l(K)\neq\0$. Is it possible to change the system $\eS$ slightly  to such system $\eS'=\{S_1',...S_m'\}$, that its attractor $K'$ satisfies the condition  $S_k'(K)\cap S_l'(K)=\0$?\\
In the above situation such statements as \cite[Theorem 8.1, Corollary 8.2] {Fal} cannot help, because the transformation from $S_k(K)$ to $S_k'(K')$ is not even a homeomorphism.\\
The following Theorem gives a way out of this situation:

\begin{thm}\label{genpos}
Let $D,L_1, L_2$  be compact  metric spaces. \\
Let \ \  $\varphi_i(\xi,x):D\times L_i\to \mathbb R^n$  be  continuous maps,  such that\\
 {\bf (a)} they  are $\alpha$-H\"older  with  respect  to $x$;
 and\\  
 {\bf (b)} there is  $M>0$ such that  for  any $x_1\in L_1, x_2\in L_2$, $\xi,\xi'\in D$\\ the  function \quad  $\Phi(\xi,x_1,x_2)=\varphi_1(\xi,x_1)-\varphi_2(\xi,x_2)$ \quad satisfies  the inequality 
\begin{equation}\label{antilip }\left\|\Phi(\xi',x_1,x_2)-\Phi(\xi,x_1,x_2)\right\|\ge M|\xi'-\xi|    \end{equation}Then  $\Delta=\{\xi\in D|\ \   \varphi_1(\xi,L_1)\cap\varphi_2(\xi,L_2)\neq\0\} $ is a compact set such that \begin{equation}\label{dimda }\dim_H \Delta \leq \min\left\{\dfrac{\dim_H L_1\times L_2}{\alpha},d\right\}\end{equation} 
\end{thm}
{\bf Proof:  }
Consider the set $\tilde \Da=\{(\xi,x_1,x_2)\in D\times L_1\times L_2:\ \Phi(\xi,x_1,x_2)=0\}$ of all  coordinate triples corresponding to intersection points of sets $\varphi_1(\xi,L_1)$ and $\varphi_2(\xi,L_2)$. 
Let $\Da$ be its projection to $D$ and $\Da_L$ be its projection to $L_1\times L_2$.\\
 The  set $\tilde \Da$ is  a compact subset of  $D\times L_1\times L_2$. \\
  
Consider  the  projections $\pi_D:\tilde \Da\to\Da\IN D$  and $\pi_L:\tilde \Da\to\Da_L\IN \ L_1\times L_2$.  The projection $\pi_L$ is a homeomorphism. Indeed, the  condition {\bf(b)} of  the Theorem implies that
$\pi_L:\tilde \Da\to\Da_L$ is  a bijection, so $\pi_L$ is a continuous bijection of  compact sets.
 
Define the map $g=\pi_D\circ\pi_L^{-1}:\Da_L\to\Da$. This map is $\al$-H\"older with  respect  to $x_1$ and $x_2$. \\To show  this, take some $\xi=g(x_1,x_2)$ and $\xi'=g(x'_1,x'_2)$. These two equalities are equivalent to
$\Phi(\xi',x'_1,x'_2)=\Phi(\xi,x_1,x_2)=0$, which  implies   the following  double inequality:
 $$C(\|x_1'-x_1\|^\alpha+\|x_2'-x_2\|^\alpha)\ge\|\Phi(\xi',x'_1,x'_2)-\Phi(\xi',x_1,x_2)\|=$$ $$=\|\Phi(\xi,x_1,x_2)-\Phi(\xi',x_1,x_2)\|\ge M|\xi'-\xi|$$
 
 Comparing the first and  the last part of this inequality we see that the function $g(x_1,x_2)$ is $\al$-H\"older with respect to both of  its  variables.
Since $\dim_H\Da_L\le \dim_H (L_1\times L_2)$ and $\Da=g(\Da_L)$, 
$$\dim_H\Delta\le\min\left\{\dfrac{\dim_H L_1\times L_2}{\alpha},d\right\}\qquad \blacksquare$$

\medskip

\subsubsection{\bf Displacement Theorem}

To evaluate the displacement  $|\pi (\sigma)-\pi' (\sigma)|$ of elements $x=\pi (\sigma)$ of the set $K_{pq}$ under the transition to the set $K_{pq'}$ we use the following Displacement Theorem, which can be considered as a modification of Barnsley Collage Theorem \cite{Barnsley}:

\medskip

\begin{thm}\label{collage}

Let  $\eS=\{S_1,...,S_m  \}$ and $\eS'=\{S'_1,...,S'_m  \}$   be  two  systems of  contractions in $\rr^n$. Let  $\pi:\ I^\infty\to K$ and  $\pi':\ I^\infty\to K'$ be the address maps with $I=\{1,...,m\}$. Suppose $V$ is  such compact set, that for  any $i=1,...,m$, $S_i(V)\IN V$  and  $ S'_i(V)\IN V$.\\
Then, for  any $\sa\in \io$, 
\begin{equation}\label{deltapi}\|\pi(\sa)-\pi'(\sa)\|\le\dfrac{\da}{1-q},\end{equation}
where   $p=\max\limits_{i\in I} (\Lip S_i,\Lip S'_i)$  and $\da=\max\limits_{x\in V, i\in I} \| S'_i(x)- S_i(x)\|$.
\end{thm}
{\bf Proof:  }\ 
Take $\sigma=i_1 i_2\dots$ and denote $\sigma_k=i_k i_{k+1}\dots$. Since $\pi(\sigma_{k})=S_{i_k}\pi(\sigma_{k+1})$,\\ 
$\| \pi(\sigma_{k})-\pi'(\sigma_{k})\| \le 
\|S_{i_k}\pi(\sigma_{k+1})-S_{i_k}\pi'(\sigma_{k+1}) \|+
\|S_{i_k}\pi'(\sigma_{k+1})-S_{i_k}'\pi'(\sigma_{k+1}) \|$, so\\
$\| \pi(\sigma_{k})-\pi'(\sigma_{k})\|\le p \| \pi(\sigma_{k+1})-\pi'(\sigma_{k+1})\|+\delta$ for any $k\in \nn$.\\
Therefore $\|\pi(\sigma)-\pi'(\sigma)\|\le p^{n+1} \|\pi(\sigma_{n+1})-\pi'(\sigma_{n+1})\|+ \delta\sum\limits_{k=0}^{n}p^k$, which becomes (\ref{deltapi}) as $k$ tends to $\8$. 
\vse \\

\medskip

{\bf Notation: The space $I^\8_a$.} \quad Let $0<a<1$ and $I^\8_a$ be the space $I^\8$ supplied with the metrics $\rho_a(\sa,\tau)=a^{s(\sa,\tau)}$, where
 $s(\sa,\tau)=\min\{k:\sa_k=\tau_k\}-1$.
 
 This metrics turns $I^\8$ to a self-similar set having Hausdorff dimension $\dim_HI^\8_a=-\dfrac{\log 4}{\log a}$.  Particularly, if  $0<a<\dfrac{1}{16}$,\  then    $\dim_HI^\8_a<1/2$.\\

 We use the space $I^\8_a$ to parametrise the sets $K_{pq}$ and further apply this in Theorem \ref{genpos}:\\

\begin{lem}\label{tech1}
Let $p,q\in(0,a)$ and $a<\dfrac{1}{16}$. Then $\pi_{pq}:I^\8_a\to K_{pq}$ is a 1-Lipschitz map.
\end{lem}
{\bf Proof:  }\ 
Take $\sa,\tau\in I^\8$ and let $s(\sa,\tau)=k$. Then $\rho_a(\sa,\tau)=a^k$.
There is $\bi\in I^k$ such that $\bi\sqsubset \sa$ and $\bi\sqsubset \tau$, so both $\pi_{pq}(\sa)$ and $\pi_{pq}(\tau)$ are contained in $S_\bj(K)$, whose diameter is equal to $\Lip(S_\bj)<a^k$. Therefore $\dfrac{|\pi_{pq}(\sa)-\pi_{pq}(\tau)|}{\rho_a(\sa,\tau)}<1$.
\vse

\begin{lem}\label{tech2}
Let $p<1/16$, $D_{mn}(p)=\left\{q<1/16:\ \dfrac{15}{16}\le\dfrac{q^n}{p^m}\le \dfrac{16}{15} \right\}$.\\
Let $\fy_1(q, \sa)=S_1^m S_i \pi_{pq}(\sa)$ and $\fy_2(q, \tau)=S_2^n S_j \pi_{pq}(\tau)$, where $i,j\in\{3,4\}$.\\
Then  for any $\sa,\tau\in I^\8$ and for any $q,q' \in D_{mn}(p)$: 
\begin{equation}\label{ge11d }|\fy_1(q, \sa)-\fy_2(q, \tau)-\fy_1(q', \sa)+\fy_2(q', \tau)|>11p^m|q'-q| \end{equation}
\end{lem}
{\bf Proof:  }\   Take  $\eS_{pq}=\{S_1,S_2,S_3,S_4\}$ and $\eS_{pq'}=\{S'_1,S'_2,S'_3,S'_4\}$, where $q,q' \in D_{mn}(p)$.

Notice that $S'_2(x)=q'x$ and $S_4'(x)=1-q'+q'x$, while $S_1'=S_1$  and $S_3'=S_3$.\\
Let $x=\pi_{pq}(\sa)$, $x'=\pi_{pq'}(\sa)$, $y=\pi_{pq}(\tau)$, $y'=\pi_{pq'}(\tau)$ be the images of $\sa,\tau$ in $K_{pq}$ and $K_{pq'}$.\\

Denote $\da=|q-q'|$, $\da_{1i}= S_1^mS'_ix'-S_1^mS_ix$ and  $\da_{2j}={S'}_2^nS'_jy'-S_2^nS_jy$.\\ 
 
Thus, we have to find a lower bound for $|\da_{2j}-\da_{1i}|$ valid for any $i,j\in\{3,4\}$.

It follows from Theorem \ref{collage} that $|x-x'|$ and $|y-y'|$ do not exceed $\dfrac{16\da}{15}$.\\
Without loss of generality we can take $q<q'$. By Lagrange theorem for any $k\in\nn$, 
\begin{equation}\label{lagr }k q^{k-1}\delta \le |q'^k-q^k|\le kq'^{k-1} \delta \end{equation}

Let us evaluate $|S'_i(x')-S_i(x)|$ for $i=3,4$.\\  $|S'_3(x')-S_3(x)|=|p(x'-x)|\le \dfrac{p\cdot16\da}{15}<\dfrac{\da}{15}$.\\
$|S'_4(x')-S_4(x)|=|q'(1-x')-q(1-x)|\le |\da(1-x')|+|q(x'-x)|<\da+\dfrac{q\cdot16\da}{15}<\dfrac{16\da}{15}$.\\

Writing $\da_{2j}=({S'_2}^nS_j'(y')-{S_2}^nS_j'(y'))+({S_2}^nS_j'(y')-{S_2}^nS_j(y))$, we notice that by (2), 
 \begin{equation*}|{S'_2}^nS_j'(y')-{S_2}^nS_j'(y')|\ge\dfrac{15}{16}nq^{n-1}\da\end{equation*} 
 
 So we have $|\da_{2j}|>\left(\dfrac{15n}{16q}-\dfrac{16}{15}\right)q^n\da$.\\ Taking $n\ge1$ and $q^n\ge\dfrac{15}{16}p^m$,  we get $|\da_{2j}|>(14n-1)p^m\da\ge 13p^m\da$.   \\
For $\da_{1i}={S_1}^mS_i'(x')-{S_1}^mS_i(x)$ we see that $|\da_{1i}|<\dfrac{16p^m\da}{15}<2p^m\da$.\\
Therefore $|\da_{2j}-\da_{1i}|>11p^m\da$. \vse \\

\medskip

\subsubsection{\bf Almost all $K_{pq}$ are twofold Cantor sets.}

\begin{thm}\label{tech3}
Let $p\in(0,1/16)$. Then for any $m,n\in\mathbb N$ the set $\Delta_{mn}(p)=\{q\in(0,1/16):\ S_1^m(A)\cap S_2^n(A)\neq\0 \}$ is closed  and nowhere dense in $(0,1/16)$.
\end{thm}
{\bf Proof:  }\ 
Take some $a\in(0,1/16)$ and $p\in(0,a)$.\\
Consider the functions $\fy_1(q,\sa)=S_1^m S_i \pi_{pq}(\sa)$ и $\fy_2(q, \sa)=S_2^n S_j \pi_{pq}(\sa)$, where $i,j\in\{3,4\}$ as  maps from $I^\8_a$ to $K_{pq}$. It follows from Lemma \ref{tech1} that they are 1-Lipschitz with respect to $\sa$, and from Lemma   \ref{tech2} it follows that if $q,q'\in D_{mn}(p)\cap(0,a)$ and $\Phi(q,\sa,\tau)=\fy_1(q,\sa)-\fy_2(q,\tau)$ then: 
 \begin{equation}\label{deltaF}\left\|\Phi(q',\sa,\tau)-\Phi(q,\sa,\tau)\right\|\ge 11p^m\|q'-q\|\end{equation}

Consider the set $D_{mn}(p)=\left\{q<1/16:\ \dfrac{15}{16}\le\dfrac{q^n}{p^m}\le \dfrac{16}{15} \right\}$. For any $a\in(0,1/16)$ the set $D_{mn}(p)\cap(0,a)$ is either closed interval, either can be covered by countable many closed intervals.
Applying Theorem \ref{genpos} to these closed intervals we get that the set $\Delta_{mn}(p)\cap(0,a)$ is closed in $(0,a)$ and its dimension is not greater than $2\dim_H I^\8_a<1$. Therefore $\Da_{mn}(p)$ is closed in $(0,1/16)$ and its dimension is less or equal to 1, so it has zero $H^1$-measure and is nowhere dense in $(0,1/16)$.\\
\vse 

\medskip

\begin{cor}\label{tech4}
The set
$\Delta_{mn}'=\{(p,q):\ p,q<1/16,\ S_1^m(A)\cap S_2^n(A)\neq\0 \}$ is a null-measure closed subset in $(0,1/16)^2$.
\end{cor}
{\bf Proof:  }
Define  a function $\Psi: (0,1/16)^2\times (I^\8)^2 \to \rr$   by $\Psi(p,q,\sa,\tau)=|S_1^m S_i \pi_{pq}(\sa)-S_2^n S_j \pi_{pq}(\tau)|$. It is continuous, therefore the set $\Psi^{-1}(\{0\})$ is closed in $(0,1/16)^2\times (I^\8)^2$.  The projection of this set to 
$(0,1/16)^2$ is $\Da'_{mn}$. Since $(I^\8)^2$ is compact, $\Da'_{mn}$ is closed in $(0,1/16)^2$.\\

By Fubini's Theorem, 2-dimensional Lebesgue measure of the set $\Delta_{mn}'$ is equal to\\ $\iint\limits_{(0,1/16)^2} \chi (p,q) dp\ dq=\int\limits_{0}^{1/16} dp \int\limits_{0}^{1/16} \chi(p,q) dq$, where $\chi(p,q)$ is a characteristic function of the set $\Delta_{mn}'$. By Theorem \ref{tech3},   $\int\limits_{0}^{1/16} \chi(p,q) dq=0$ for any $p\in (0,1/16)$. Therefore the set $\Delta_{mn}'$ has zero measure in $(0,1/16)^2$. \vse

\begin{thm}\label{fullmes}
The set $\eK$ of those $(p,q)\in \eV=(0,1/16)^2$, for which $K_{pq}$ is a twofold Cantor set, has full measure in $\eV$, and its complement is uncountable and dense in  $\eV$.
\end{thm}
{\bf Proof:  }\ 
By Corollary \ref{tech4},   $\Delta_{mn}'$ has  zero Lebesgue 2-measure in $\eV$ for any $m,n\in\mathbb N$. Since $\Delta=\bigcup\limits_{m,n=1}^{\infty} \Delta_{mn}'$, the same is true for $\Da$ and its complement  $\mathcal{K}$ has full measure in $\eV$.\\
From the other side, the set $\Delta_0=\{(p,q)\in (0,1):\ \dfrac{\log p}{\log q} \in\mathbb Q\}\cap \eV \subseteq \Delta$ is a countable union of curves $q=p^s$ in $\eV$.
 The set $\{p^s:\ s\in\mathbb Q_{+}\}$ is dense in $(0,1/16)$ for any $p\in(0,1/16)$, therefore $\Delta$ is uncountable and dense in $\eV$.
\vse

\section*{Acknowledgments.}

The authors would like to thank Christoph Bandt for very intense long-lasting discussions at all stages of writing the paper.

The research of authors was supported by Russian Foundation of Basic Research grants No.~16-01-00414, No.~18-501-51021.

\end{document}